\documentclass[10pt]{article}

\usepackage{amsmath}
\usepackage{amsfonts}
\usepackage{amssymb}

\usepackage{color}

\newtheorem{lemma}{Lemma}
\newtheorem{theorem}{Theorem}

\newtheorem{corollary}{Corollary}

\title{Any decreasing cycle--convergence curve is possible for restarted GMRES}
\author{Eugene Vecharynski and Julien Langou}

\begin{document}

\maketitle

\begin{abstract}
Given a matrix order $n$, a restart parameter $m$ ($m < n$),
a decreasing positive sequence $f(0) > f(1) > \ldots > f( q ) \geq 0$, where $q < n/m$,
it is shown that there exits an $n$-by-$n$ matrix $A$ and a vector $r_0$ with
$\|r_0\|=f(0)$ such that $\|r_k\|=f(k)$, $k=1,\ldots,q$, where $r_k$ is the
residual at cycle $k$ of restarted GMRES with restart parameter $m$ applied to the linear
system $Ax=b$, with initial residual $r_0=b-Ax_0$. Moreover, the matrix $A$ can
be chosen to have any desired eigenvalues. We can also construct arbitrary
cases of stagnation; namely, when $f(0) > f(1) > \ldots > f(i) = f(i+1) \geq 0  $
for any $i <q $. The restart parameter can be fixed or variable.
\end{abstract}

\section{Introduction}\label{sec:intro}

We consider the \textit{generalized minimal residual method}
(GMRES)~\cite{SaadSchultz:86} for solution of a nonsingular non-Hermitian
systems of linear equations

\begin{equation}\label{eqn:sys}
A x = b, \quad A \in \mathbb{C}^{n \times n}, \quad b \in \mathbb{C}^n.
\end{equation}

For a few class of matrices, some convergence estimates are available for
\textit{restarted} GMRES and \textit{full} GMRES.  For example for real
positive definite matrices (that is, for matrices $A$ for which $H = (A+A^H)/2$
is symmetric positive definite, or, equivalently, for matrices $A$ for which
$x^H A x > 0$ for any nonzero $x \in \mathbb{R}^n$), the Elman's
bound~\cite{Eisenstat.Elman.Schultz.83,Elman.82,Greenbaum:97,SaadSchultz:86}
can be stated as follows
$$\| r_k \|^2 \leq (1-\rho)^k \| r_0 \|^2 \quad \mbox{where } 0< \rho \equiv (\lambda_{min}(H)/\|A\|)^2\leq 1. $$
The latter guarantees linear convergence of GMRES($m$) for any value of $m\geq
1$ for a positive definite matrix. Improvements and generalizations of this
bound can be found
in~\cite{Beckermann.Goreinov.Tyrtyshnikov.05,SimonciniSzyld:08,Zitko:00}.

For normal matrices the convergence of both full and restarted GMRES is well studied. 
In particular, the convergence of \textit{full} GMRES for normal matrices is known to be linear
and there exist convergence estimates governed solely by the spectrum of
$A$~\cite{Simoncini.Szyld.05a,vanderVorstVuik:93}.
The convergence of \textit{restarted} GMRES for normal matrices, however, is
sublinear~\cite{BakerJessupKolev:09,VechLangou:09}.  
The current paper is concerned with the general case.

For the general case, the following theorem proves that we can not prove convergence results
based on the spectrum of the coefficient matrix alone.
\begin{theorem}\emph{(Greenbaum, Pt\'{a}k, and Strako\v{s}, 1996,~\cite{GreenbaumPtakStrakos:96})}
\label{theo:GreenbaumPtakStrakos:96}
Given a nonincreasing positive sequence $f(0) \geq f(1) \geq \cdots \geq f(n-1)
> 0$, there exists an $n$-by-$n$ matrix $A$ and a vector $r_0$
with $\| r_0 \| = f(0)$ such that $f(k) = \| r_k \|$, $k=1, \ldots , n-1$,
where $r_k$ is the residual at step $k$ of the GMRES algorithm applied to the
linear system $Ax=b$, with initial residual $r_0 = b - A x_0$. Moreover, the
matrix $A$ can be chosen to have any desired eigenvalues.
\end{theorem}

This result states that, in general, eigenvalues alone do not determine the convergence of
\textit{full} GMRES.  Assuming that the coefficient matrix $A$
is diagonalizable, some characterizations of the convergence of full GMRES rely
on the condition number of the eigenbasis~\cite{vanderVorstVuik:93}.  Other
characterizations of the convergence of full GMRES rely on
pseudospectra~\cite{NachtigalReddyTrefethen:92}.  More commonly, the field of
values is used~\cite{Beckermann.Goreinov.Tyrtyshnikov.05,Eisenstat.Elman.Schultz.83,Elman.82,Greenbaum:97,SaadSchultz:86,SimonciniSzyld:08,Zitko:00}.
A discussion on how descriptive some of these bounds are is given 
by Embree~\cite{Embree.99}.
 
The main result of this paper is given in the abstract. We will repeat it here
\begin{theorem}
\label{theo:main}
Given a matrix order $n$, a restart parameter $m$ ($m < n$),
a decreasing positive sequence $f(0) > f(1) > \ldots > f( q ) \geq 0$, where $q < n/m$,
there exits an $n$-by-$n$ matrix $A$ and a vector $r_0$ with
$\|r_0\|=f(0)$ such that $\|r_k\|=f(k)$, $k=1,\ldots,q$, where $r_k$ is the
residual at cycle $k$ of restarted GMRES with restart parameter $m$ applied to the linear
system $Ax=b$, with initial residual $r_0=b-Ax_0$. Moreover, the matrix $A$ can be chosen
to have any desired eigenvalues.
\end{theorem}
Section~\ref{sec:constr} contains a proof of Theorem~\ref{theo:main}.
Theorem~\ref{theo:main} is to \textit{restarted} GMRES what Theorem~\ref{theo:GreenbaumPtakStrakos:96} is to \textit{full} GMRES. 
The proof we provide is constructive and directly inspired by the article of 
Greenbaum, Pt\'{a}k, and Strako\v{s}~\cite{GreenbaumPtakStrakos:96}.
Although Greenbaum, Pt\'{a}k, and Strako\v{s} laid the path, there  
are several specific difficulties ahead in the case of \textit{restarted} GMRES.

\textit{Full} GMRES has a nonincreasing convergence (for any $i\geq 0$, $f(i)\geq
f(i+1)$) and it computes the exact solution in at most $n$ steps ($f(n)=0$).  It is
remarkable that Greenbaum, Pt\'{a}k, and Strako\v{s} are able to characterize any
{\em admissible} convergence for GMRES.  (See assumptions on $f$ in
Theorem~\ref{theo:GreenbaumPtakStrakos:96}.)
At the same time we would like to note that the cycle--convergence of \textit{restarted} GMRES 
can have two {\em admissible} scenarios: either for any $i$, $f(i) > f(i+1)$, in other words, the
cycle--convergence is decreasing; or there exits $s$ such that $f(i) > f(i+1)$ for any $i<s$, 
and then for any $i > s,$ $f(i)=f(s)$, in other words, if restarted GMRES
stagnates at cycle $s+1$, it stagnates forever. Theorem~\ref{theo:main}
considers the first case (decreasing cycle--convergence).  In
Section~\ref{sec:stagnation}, we consider the second case (stagnation). Therefore
with Theorem~\ref{theo:main} and Section~\ref{sec:stagnation}, we prove that any
{\em admissible} cycle--convergence curve is possible for the $q$ first
cycles of \textit{restarted} GMRES.

As mentioned above, the maximum number of iterations of \textit{full} GMRES is at most $n$, 
and the method delivers the exact solution in a finite number of steps. 
\textit{Restarted} GMRES, however, may never provide the exact solution.
It will (hopefully) decrease the residual norm at each cycle, that is, provide a more and more accurate 
approximation to the exact solution. With $n^2$ parameters in $A$
and $n$ parameters in $b$ we are not able to control the convergence for an
infinite amount of cycles. For this reason, it is natural to consider only
the first $q < n/m$ initial GMRES($m$) cycles. Actually, we provide the same level of control as
Greenbaum, Pt\'{a}k, and Strako\v{s}: $n$ iterations (or $q$ cycles with
$q < n/m$) and $n$ eigenvalues. 

In Section~\ref{sec:variable_restart}, we generalize the result
given by Theorem~\ref{theo:main} and Section~\ref{sec:stagnation} for the case 
of variable restart parameters. The sequence of restart parameters $m_k$ needs to be known a priori. 
We show that GMRES($m_k$) can produce any admissible cycle--convergence curve   
at the $q$ initial cycles, regardless of the spectrum of the coefficient matrix,
where $q$ is such that $\sum_{i=1}^q m_k <n$. We note that our construction
can be a reasonable tool for generating examples/counter-examples for
different strategies for varying the restart parameter, 
e.g.~\cite{BakerJessupKolev:09}. 

The cycle--convergence of restarted GMRES for normal matrices is
sublinear~\cite{BakerJessupKolev:09,VechLangou:09}. However, for general
matrices, through Theorem~\ref{theo:main}, one can expect any convergence curve. 
In particular, it is possible to construct matrices for which the convergence of 
GMRES($m$) is fast (e.g.  superlinear). This
relates to the observations of Zhong and Morgan~\cite{ZhongMorgan:08} who
report superlinear cycle--convergence for their particular cases 
of nonnormal matrices, as well as to~\cite{VechLangou:09}, where it is shown that
the cycle--convergence can become superlinear as the coefficient matrix departs from normality.

In a pedagogical paper, Embree~\cite{Embree:03} presents a 3-by-3
linear system of equations and attempts to solve it with GMRES(1) 
and GMRES(2). While GMRES(1) converges to the exact solution in 
3 cycles, GMRES(2) (almost) stagnates. Our main result, basically,
reaffirms this intuition in the sense that the increase in the restart parameter
(and thus, in the computational complexity at each cycle) does not
necessarily imply a faster convergence.

In order to improve the convergence of restarted GMRES, several
techniques~\cite{Baglama.Calvetti.Golub.Reichel:99,Chapman.Saad:1997,Duff.Giraud.Langou.Martin:05,Erhel.Burrage.Pohl:96,Morgan.02}
have been proposed which consist of augmenting (or enriching) the Krylov space
with eigenvectors or, alternatively, deflating some of the eigenvalues from the
spectrum of the original matrix $A$. The eigenvalues targetted are the ones the
closest from zero. These techniques have proved effective and the convergence
of restarted GMRES is, in practice, greatly improved.  Theorem~\ref{theo:main}
states that, in the general case, eigenvalues alone do not determine the
convergence of restarted GMRES, therefore it is hard to provide a theorical
justification for the choice of removing the eigenvalues the closest from zero.
A beginning of theoretical understanding has been provided by
Z{\'{\i}}tko~\cite{Zitko:08}.

We have generated two Matlab functions that correspond to
Theorem~\ref{theo:GreenbaumPtakStrakos:96} and Theorem~\ref{theo:main}.
Given a matrix size, a restart parameter, a convergence curve and a spectrum, 
we construct the appropriate matrix and right-hand side.
See:
\texttt{http://www-math.cudenver.edu/\string ~eugenev/edf.software/anycurve/}.

The main message that we would like our readers to retain from this paper is
that in the context of GMRES($m$), for a certain number of initial cycles,
any convergence curve is possible independently of the spectrum of the coefficient
matrix.  This means that eigenvalues alone do not determine the convergence of
\textit{restarted} GMRES.

\section{Constructive proof of Theorem~\ref{theo:main}}\label{sec:constr}

Let $n$ be a matrix order and $m$ a restart parameter ($m < n$),
$\Lambda = \lbrace \lambda_1, \lambda_2, \dots \lambda_n \rbrace \subset \mathbb{C} \setminus \lbrace 0 \rbrace$ 
be a set of $n$ nonzero complex numbers,   and $\{ f(k) \}_{k=0}^q$
be a decreasing sequence of positive real numbers, $q<n/m$.


In this section we construct a matrix $A \in \mathbb{C}^{n \times n}$
and an initial residual vector $r_0 = b - Ax_0 \in \mathbb{C}^{n}$
such that GMRES($m$) applied to the system (\ref{eqn:sys}) with the initial
approximate solution $x_0$,
produces a sequence $\{x_k\}_{k=1}^q$ of approximate solutions
with corresponding residual vectors $\{r_k\}_{k=0}^q$ having the prescribed norms: 
$\| r_k \| = f(k)$. Moreover the spectrum of $A$ is $\Lambda$.

\subsection{Outline of the proof}\label{subsec:gen_schem}

The general approach described in this paper is similar to the approach of
Greenbaum, Pt\'{a}k, and Strako\v{s}~\cite{GreenbaumPtakStrakos:96}: we
fix an initial residual vector,
construct an appropriate basis of $\mathbb{C}^n$ and use this basis to define a
linear operator $\mathcal{A}$. This operator is represented by the matrix $A$
in the canonical basis. It has the prescribed spectrum and provides the desired
cycle--convergence at the first $q$ cycles of GMRES($m$).  However, the
presence of restarts somewhat complicates the construction: the choice of the
basis vectors, as well as the structure of the resulting operator
$\mathcal{A}$, becomes less transparent.  Below we describe our
three-step construction.

At the \textit{first step} we construct $q$ sets of vectors 
$\mathcal{W}_m^{(k)} = \lbrace w_1^{(k)}, \dots ,w_m^{(k)} \rbrace$, 
$k = 1,\smallskip \dots \smallskip, q$, each set $\mathcal{W}_m^{(k)}$ is the orthonormal basis of 
the Krylov residual subspace $A\mathcal{K}_m\left( A,r_{k-1} \right)$ generated at 
the $k$-th GMRES($m$) cycle such that
\begin{equation}\label{eqn:span}
\mbox{span}~\mathcal{W}_j^{(k)} = A \mathcal{K}_j\left( A,r_{k-1} \right),    
\quad j = 1,  \dots , m.
\end{equation}
(With this definition, $\mathcal{W}_m^{(k)}$ is defined up to multiplication by a complex number of unit modulus.)

The orthonormal basis $\mathcal{W}_m^{(k)}$ needs to be chosen in order to generate residual vectors $r_k$ 
with the prescribed norms $f(k)$ at the end of each
cycle subject to the additional requirement that the set of $mq+1 (\leq n)$ vectors
\begin{equation}\label{eqn:indep}
\mathcal{\overline{S}} = \lbrace r_0, w_1^{(1)}, \dots , w_{m-1}^{(1)}, r_1, w_1^{(2)},  \dots , w_{m-1}^{(2)},
\dots , r_{q-1}, w_1^{(q)},  \dots , w_{m-1}^{(q)}, r_q \rbrace 
\end{equation}
is linearly independent.

Once we have the set $\mathcal{\overline{S}}$, we will complete it to have a
basis for $\mathbb{C}^n$.  When the number of vectors in
$\mathcal{\overline{S}}$ is less than $n$, a basis $\mathcal{S}$ of
$\mathbb{C}^n$ is obtained by completion of $\mathcal{\overline{S}}$ with a set
$\mathcal{\widehat{S}}$ of $n - mq - 1$ vectors, i.e.  $\mathcal{S} = \lbrace
\mathcal{\overline{S}}, \mathcal{\widehat{S}}\rbrace$.  This will provide a
representation of $\mathbb{C}^n$ as the direct sum
\begin{equation}\label{eqn:direct_sum}
\mathbb{C}^n = \mbox{span}~\mathcal{S} =
\mbox{span}\lbrace r_0, \mathcal{W}_{m-1}^{(1)} \rbrace \oplus  \dots \oplus \mbox{span}\lbrace r_{q-1}, \mathcal{W}_{m-1}^{(q)} \rbrace \oplus \mbox{span}\lbrace r_q, \mathcal{\widehat{S}} \rbrace.
\end{equation}
The latter translates in terms of Krylov subspaces into
\begin{equation*}
\mathbb{C}^n = \mbox{span}~\mathcal{S} = 
\mathcal{K}_m\left( A,r_0 \right) \oplus \dots \oplus \mathcal{K}_m\left( A,r_{q-1} \right) \oplus \mbox{span}\lbrace r_q, \mathcal{\widehat{S}} \rbrace. 
\end{equation*}

At the \textit{second step} of our construction, we define a linear operator 
$\mathcal{A}: \mathbb{C}^n \longrightarrow \mathbb{C}^n$ with spectrum $\Lambda$
which generates the Krylov residual subspaces in Eq.~(\ref{eqn:span}) at each GMRES($m$) cycle, 
by its action on the basis vectors $\mathcal{S}$, such that the desired matrix $A$ 
is the operator $\mathcal{A}$'s representation in the canonical basis. 
The \textit{third step} accomplishes the construction by a similarity transformation.

The two following subsections are concerned with the question if 
(\ref{eqn:span})--(\ref{eqn:indep})--(\ref{eqn:direct_sum}) and the definition
of the operator $\mathcal{A}$ with the prescribed spectrum is
actually possible. 

\subsection{Step 1: Construction of a sequence of Krylov subspaces which provide the prescribed cycle--convergence}\label{subsec:subsp} 

At the $k$-th GMRES($m$) cycle, the residual vector $r_k$ satisfies the
following minimality condition:

\begin{equation}\label{eqn:gmres} 
\| r_k \| = \displaystyle\min_{u \in A \mathcal{K}_m (A,r_{k-1})} \| r_{k-1} - u\|.          
\end{equation}

We assume that each set $\mathcal{W}_m^{(k)}$ is
an orthonormal basis of a corresponding Krylov residual subspace $A \mathcal{K}_m\left( A,r_{k-1} \right)$,
therefore
the condition~(\ref{eqn:gmres}) implies

\begin{equation}\label{eqn:new_res} 
r_k = r_{k-1} - \sum_{j=1}^m  \langle r_{k-1},w_j^{(k)} \rangle w_j^{(k)}, \quad k = 1,  \dots , q.
\end{equation}

At this stage, in order to simplify the forthcoming justification of the linear independence of the set
$\mathcal{\overline{S}}$,
we impose a stricter requirement on the residual change inside the cycle. We will require that
the residual vector $r_{k-1}$ remains constant
during the first $m-1$ inner steps of GMRES and is reduced only at the last, $m$-th,
step. Thus, the equality in (\ref{eqn:new_res}) can be written as
\begin{equation}\label{eqn:simple_new_res} 
 r_{k} = r_{k-1} - \langle r_{k-1},w_m^{(k)} \rangle w_m^{(k)},       
\quad k = 1,  \dots , q.
\end{equation}
This implies that the vectors $w_j^{(k)}$, $j = 1,  \dots, m-1$,
are orthogonal to the residual vector $r_{k-1}$, i.e.
\begin{equation}\label{eqn:orth} 
\langle r_{k-1},w_j^{(k)}\rangle = 0,
\quad j = 1,  \dots , m-1,       
\quad k = 1, \dots , q.
\end{equation}

From Eq.~(\ref{eqn:simple_new_res}), using the fact that $r_k \perp w_m^{(k)}$ and the Pythagorean theorem, 
we obtain
\begin{equation*}
\vert \langle r_{k-1},w_m^{(k)} \rangle \vert =  \sqrt{\Vert r_{k-1} \Vert^2 - \Vert r_k \Vert^2},
\quad k = 1,  \dots , q.
\end{equation*}
We rewrite the expression above in terms of cosines of angles $\psi_k = \angle(r_{k-1},w_m^{(k)})$ by
prescribing the expected values $f(k)$ for the norms of the residuals. We get
\begin{equation}\label{eqn:cos} 
\cos \psi_{k}  =  \frac{\sqrt{ f(k-1)^2 - f(k)^2}}{f(k-1)} \in (0,1),       
\quad k = 1, \dots , q.
\end{equation}

This latter equation means that, if we are given $r_{k-1}$, one way to ensure
the desired cycle--convergence at cycle $k$ of GMRES($m$) is to choose the unit
vectors $w_j^{(k)}$ such that (\ref{eqn:simple_new_res})--(\ref{eqn:cos}) holds.


In the following lemma, we show constructively that the described approach
(\ref{eqn:simple_new_res})--(\ref{eqn:cos}) leads to an appropriate
set $\mathcal{\overline{S}}$.

 
\begin{lemma}\label{lemm:1}
Given an initial vector $r_0$, $\Vert r_0 \Vert = f(0)$, 
there exist vectors $r_k$, $\Vert r_k \Vert = f(k)$ and orthonormal sets $\mathcal{W}_m^{(k)}$
such that Eq. (\ref{eqn:simple_new_res}), (\ref{eqn:orth}) and (\ref{eqn:cos}) hold, and 
the set
$\mathcal{\overline{S}}$
 is linearly independent, $k = 1, \dots, q < n/m$.
\end{lemma}
\textbf{Proof.} The proof is by induction. 

Let $k = 1$. Given the initial vector $r_0$, $\Vert r_0 \Vert = f(0)$, we pick
$\mathcal{W}_{m-1}^{(1)} = \lbrace w_1^{(1)},  \dots , w_{m-1}^{(1)} \rbrace$ an 
orthonormal set in $r_0^{\perp}$ in order to satisfy Eq.~(\ref{eqn:orth}).
The set $\lbrace r_0, \mathcal{W}_{m-1}^{(1)} \rbrace$ is linearly independent.

In order to choose the unit vector $w_m^{(1)}$ 
orthogonal to the previously constructed vectors $\mathcal{W}_{m-1}^{(1)}$ and satisfying Eq.~(\ref{eqn:cos}), 
we introduce a unit vector $y^{(1)} \in \lbrace r_0, \mathcal{W}_{m-1}^{(1)} \rbrace^{\perp}$, so that
\begin{equation*}
 w_m^{(1)} = \frac{r_0}{f(0)} \mbox{cos}\psi_1  + y^{(1)}\mbox{sin}\psi_1.
\end{equation*}

We find the vector $r_1$ by satisfying Eq.~(\ref{eqn:simple_new_res}).
Eq.~(\ref{eqn:cos}) guarantees that $\Vert r_1 \Vert = f(1)$, as desired.
Finally, we append the constructed vector $r_1$ to $\lbrace r_0, \mathcal{W}_{m-1}^{(1)}\rbrace$
and get the set $\lbrace r_0, \mathcal{W}_{m-1}^{(1)}, r_1\rbrace$, which is 
linearly independent, since, by construction, $r_1$ is not in span~$\lbrace r_0, \mathcal{W}_{m-1}^{(1)} \rbrace$.

The induction assumption is that we have constructed $k-1$ vectors $r_1, \ldots, r_{k-1}$ with the prescribed norms
$f(1), \ldots, f(k-1)$ and orthonormal sets $\mathcal{W}_{m}^{(1)}, \ldots, \mathcal{W}_{m}^{(k-1)}$,
such that the equalities (\ref{eqn:simple_new_res}), (\ref{eqn:orth}) and (\ref{eqn:cos}) hold, and 
the set 
\begin{equation}\label{eqn:intermediate}
\lbrace r_0, \mathcal{W}_{m-1}^{(1)},  \ldots , r_{k-2}, \mathcal{W}_{m-1}^{(k-1)} , r_{k-1} \rbrace                
\end{equation}
is linearly independent. We want to show that we can construct the next vector $r_k$,  
$\Vert r_k \Vert = f(k)$, and the orthonormal set
$\mathcal{W}_m^{(k)}$, satisfying Eq.~(\ref{eqn:simple_new_res}),~(\ref{eqn:orth}) and~(\ref{eqn:cos}), 
such that 
\begin{equation}\label{eqn:intermediate_next}
\lbrace r_0, \mathcal{W}_{m-1}^{(1)},  \ldots , r_{k-2}, \mathcal{W}_{m-1}^{(k-1)} , r_{k-1}, \mathcal{W}_{m-1}^{(k)}, r_k \rbrace  
\end{equation}
is linearly independent, $k \leq q$.

We start by constructing orthonormal vectors 
$\mathcal{W}_{m-1}^{(k)} = \lbrace w_1^{(k)}, \dots, w_{m-1}^{(k)} \rbrace$,
satisfying Eq.~(\ref{eqn:orth}), with the additional requirement that the set $\mathcal{W}_{m-1}^{(k)}$ is not
in the span of the previously constructed vectors given in the set~(\ref{eqn:intermediate}).
From these considerations we choose $\mathcal{W}_{m-1}^{(k)}$ as an orthonormal set 
in the orthogonal complement of (\ref{eqn:intermediate}), i.e.
\begin{equation*}
w_j^{(k)} \in \lbrace r_0, \mathcal{W}_{m-1}^{(1)},  \ldots , r_{k-2}, \mathcal{W}_{m-1}^{(k-1)} , r_{k-1} \rbrace^\perp,
\quad j = 1, \dots ,m-1.
\end{equation*}
Appending $\mathcal{W}_{m-1}^{(k)}$ to the set~(\ref{eqn:intermediate}) will give a linearly independent set.  

To finish the proof, we need to construct the vector $w_m^{(k)}$, satisfying Eq.~(\ref{eqn:cos}) and orthogonal
to $\mathcal{W}_{m-1}^{(k)}$.
For this reason we introduce a unit vector $y^{(k)}$,
\begin{equation*}
y^{(k)} \in \lbrace r_0, \mathcal{W}_{m-1}^{(1)},  \ldots , r_{k-2}, \mathcal{W}_{m-1}^{(k-1)} , r_{k-1}, \mathcal{W}_{m-1}^{(k)} \rbrace^{\perp} ,
\end{equation*}
so that $w_m^{(k)}$
\begin{equation*}
w_m^{(k)} = \frac{r_{k-1}}{f(k-1)} \mbox{cos} \psi_{k}  + y^{(k)}\mbox{sin}\psi_{k}.
\end{equation*}

We define the vector $r_k$ with Eq.~(\ref{eqn:simple_new_res}). Eq.~(\ref{eqn:cos}) guarantees 
$\Vert r_k \Vert = f(k)$. 
The set~(\ref{eqn:intermediate_next}) is linearly independent, since, by construction,
the vector $r_k$ is not in span~$\lbrace r_0, \mathcal{W}_{m-1}^{(1)},  \ldots , r_{k-2}, \mathcal{W}_{m-1}^{(k-1)} , r_{k-1}, \mathcal{W}_{m-1}^{(k)} \rbrace$.
\begin{flushright}
\hfill $\square$
\end{flushright}

\subsection{Step 2: Definition of a linear operator with any prescribed spectrum} 

So far we have shown that, given an initial residual vector $r_0$, $\|r_0\| = f(0)$, it is possible
to construct vectors $r_k$, $\|r_k\| = f(k)$, and orthonormal vectors $\mathcal{W}_m^{(k)}$, $k = 1, \ldots, q$,
satisfying Eq.~(\ref{eqn:simple_new_res}), (\ref{eqn:orth}) and (\ref{eqn:cos}), such that the set $\overline{\mathcal{S}}$ of $mq + 1$ vectors is linearly independent.


In order to define a unique linear operator, we need to have a valid basis of $\mathbb{C}^n$ on hand.
Thus, we expand the set $\overline{\mathcal{S}}$ by linearly independent vectors 
$\widehat{\mathcal{S}} = \lbrace \widehat{s}_1,  \dots , \widehat{s}_t \rbrace$, $t = n - mq - 1$:
\begin{equation}\label{eqn:basis}
\mathcal{S} = \lbrace r_0, \mathcal{W}_{m-1}^{(1)}, \dots, r_{q-1}, \mathcal{W}_{m-1}^{(q)}, r_q, \widehat{s}_1, \dots, \widehat{s}_t \rbrace, 
\end{equation}      
so that $\mathcal{S}$ is a basis of $\mathbb{C}^n$.
 
Before we define a linear operator $\mathcal{A}$, let us consider 
the set $\Lambda = \{ \lambda_1, \lambda_2, \smallskip \dots \smallskip, \lambda_n \}$ 
of nonzero numbers in the complex plane that will define $\mathcal{A}$'s spectrum. 
We split $\Lambda$ into $q+1$ disjoint subsets
\begin{equation*}
 \Lambda = \lbrace \Lambda_1, \Lambda_2,  \dots , \Lambda_q, \Lambda_{q+1} \rbrace,
\end{equation*}
such that each $\Lambda_k$, $k = 1,  \dots , q$, contains $m$ elements of $\Lambda$, and
the remaining $n - mq$ elements are included into $\Lambda_{q+1}$.

For each set $\Lambda_k$ we define a monic polynomial $p_k(x)$, such that the roots of this polynomial are exactly the elements
of the corresponding $\Lambda_k$:
\begin{eqnarray}
\label{eqn:charact_m}   p_k(x)     & = & x^{m} - \sum_{j=0}^{m-1}\alpha_j^{(k)} x^{j}, \quad k = 1, \dots ,q;       \\
\label{eqn:charact_t}   p_{q+1}(x) & = & x^{t+1} - \sum_{j=0}^{t}\alpha_j^{(q+1)} x^{j}, \quad t = n-mq-1;  
\end{eqnarray}
with $\alpha_j^{(k)}$'s being the coefficients of the respective polynomials, $\alpha_0^{(k)} \neq 0$, $k = 1,  \dots  q+1$. 
$p_k(x)$ can be considered as the characteristic polynomial of an $m$-by-$m$ matrix 
with spectrum $\Lambda_k$.
$p_{q+1}(x)$ can be considered as the characteristic polynomial of a $(t+1)$-by-$(t+1)$ matrix 
with spectrum $\Lambda_{q+1}$. 

We define the operator $\mathcal{A}: \mathbb{C}^n \longrightarrow \mathbb{C}^n$ as follow: 
\begin{eqnarray}
\nonumber   \mathcal{A} r_{k-1}    & = &     w_1^{(k)},       \\
\nonumber   \mathcal{A} w_1^{(k)}  & = &     w_2^{(k)},       \\
\nonumber                               &\vdots&                         \\
\nonumber   \mathcal{A} w_{m-2}^{(k)}   & = &  w_{m-1}^{(k)},                                       \\
\nonumber   \mathcal{A} w_{m-1}^{(k)}  & = &  -\alpha_0^{(k)} r_k + \alpha_0^{(k)} r_{k-1} + 
                                    \alpha_1^{(k)} w_1^{(k)} +  \dots +\alpha_{m-1}^{(k)} w_{m-1}^{(k)}, 
                                    \quad k = 1,\dots q;                                       \\
\nonumber  \mathcal{A} r_{q}  & = &  \widehat{s}_1,                                               \\ 
\nonumber  \mathcal{A} \widehat{s}_1  & = & \widehat{s}_2,                                       \\
\label{eqn:operator}                   &\vdots&                         																		\\                  \nonumber  \mathcal{A} \widehat{s}_{t-1}  & = &  \widehat{s}_t,                                   \\
\nonumber  \mathcal{A} \widehat{s}_t  & = &  \alpha_0^{(q+1)} r_{q} + \alpha_1^{(q+1)} \widehat{s}_{1} +  \dots 
+\alpha_{t}^{(q+1)} \widehat{s}_{t},                                             
\end{eqnarray}
where $\alpha_j^{(k)}$'s are the coefficients of polynomials (\ref{eqn:charact_m}) and (\ref{eqn:charact_t}).

The following lemma shows that, given vectors $r_k$ and orthonormal sets $\mathcal{W}_m^{(k)}$ constructed
according to Lemma \ref{lemm:1}, the linear operator $\mathcal{A}$, defined by (\ref{eqn:operator}) and represented
by a matrix $A$ in the canonical basis,
generates the desired Krylov residual subspaces given in Eq.~(\ref{eqn:span});
and the spectrum of $\mathcal{A}$ can be arbitrarily chosen.   
\begin{lemma}\label{lemm:2}
Let the initial residual vector $r_0$, $\|r_0 \| = f(0)$, as well as the  
residual vectors $r_k$ and orthonormal sets $\mathcal{W}_m^{(k)}$ be constructed according
to Lemma~\ref{lemm:1}. Let $\mathcal{S}$ be the basis of $\mathbb{C}^n$ as defined by Eq.~(\ref{eqn:basis}).
We assume a matrix $A$
to be
 the representation 
in the canonical basis
of the linear operator $\mathcal{A}$ defined by Eq.~(\ref{eqn:operator}).
Then the linear operator $\mathcal{A}$
generates the Krylov residual subspaces given in Eq.~(\ref{eqn:span}).
Moreover, $\mathcal{A}$ has the prescribed spectrum $\Lambda$.
\end{lemma}
\textbf{Proof.} Directly from the definition (\ref{eqn:operator}) of the linear operator $\mathcal{A}$,
for $k = 1, \ldots, q$, we have:
\begin{equation*}
  \mbox{span}\{\mathcal{A}r_{k-1}, \ldots, \mathcal{A}^j r_{k-1}\} = \mbox{span}~\mathcal{W}_j^{(k)},
  \quad j = 1, \dots , m-1.
\end{equation*} 
To see that, for each $k$, 
\begin{equation*}
  \mbox{span}\{\mathcal{A}r_{k-1}, \ldots, \mathcal{A}^m r_{k-1}\} = \mbox{span}~\mathcal{W}_m^{(k)},
\end{equation*} 
notice that, by Eq.~(\ref{eqn:simple_new_res}), $(- \alpha_0^{(k)} r_k + \alpha_0^{(k)} r_{k-1}) \in \mbox{span}\{w_m^{(k)}\}$. 
Thus, given the representation $A$ of the operator $\mathcal{A}$ in the canonical basis, Eq. (\ref{eqn:span})
holds for each $k$, $k = 1, \ldots, q$.

To prove that the arbitrarily chosen set $\Lambda$ is the spectrum of $\mathcal{A}$, 
let us consider the matrix $\left[ \mathcal{A} \right]_\mathcal{S} $
of the operator $\mathcal{A}$ in the basis $\mathcal{S}$:
\begin{equation}\label{eqn:matrix}
\left[ \mathcal{A} \right]_\mathcal{S} =
\left[
\begin{array}{cccccccccccccc}
 0      & 0      & \cdots & \alpha_0^{(1)}            \\
 1      & 0      & \cdots & \alpha_1^{(1)}            \\
 0      & 1      & \cdots & \alpha_2^{(1)}            \\
 \vdots & \vdots & \ddots & \vdots                    \\
 0      & 0      & \cdots & \alpha_{m-1}^{(1)}        \\
        &        &        & -\alpha_{0}^{(1)} & 0      & 0      & \cdots & \alpha_{0}^{(2)}    & & & \mathbf{0} \\
        &        &        &                   & 1      & 0      & \cdots & \alpha_{1}^{(2)}   \\
        &        &        &                   & 0      & 1      & \cdots & \alpha_{2}^{(2)}   \\
        &        &        &                   & \vdots & \vdots & \ddots & \vdots             \\
        &        &        &                   & 0      & 0      & \cdots & \alpha_{m-1}^{(2)} \\
        &        &        &                   &        &        &        & -\alpha_0^{(2)} & \ddots    \\
        &        &        &                   &        &        &        &                & \ddots & \ddots \\ 
        &        &        &\mathbf{0}   &        &       &        &    &    & -\alpha_0^{(q)} & 0 & 0 & \cdots & \alpha_0^{(q+1)} \\
        &        &        &    &        &        &        &    &    &  & 1 & 0 & \cdots & \alpha_1^{(q+1)} \\
        &        &        &    &        &        &        &    &    &  & 0 & 1 & \cdots & \alpha_2^{(q+1)} \\
        &        &        &    &        &        &        &    &    &  & \vdots & \vdots & \ddots & \vdots \\
        &        &        &    &        &        &        &    &    &  & 0 & 0 & \cdots & \alpha_t^{(q+1)} \\
\end{array}
\right]. 
\end{equation}

The matrix $\left[ \mathcal{A} \right]_\mathcal{S} $ has a block lower
triangular structure, hence $\left[ \mathcal{A} \right]_\mathcal{S} $'s
spectrum is the union of the eigenvalues of all diagonal blocks, which are the
companion
matrices corresponding to the sets $\Lambda_k$ with characteristic polynomials defined in (\ref{eqn:charact_m}) and
(\ref{eqn:charact_t}). Thus, the spectrum of $\mathcal{A}$ is $\Lambda$.
\begin{flushright}
 $\square$
\end{flushright}

\subsection{Step 3: Conclusion of the proof of Theorem~\ref{theo:main}}

Finally, we define $A$ as the representation of the operator $\mathcal{A}$ in the canonical basis: 
$\lbrace e_1, e_2, \smallskip \dots \smallskip, e_n\rbrace$,
\begin{equation}\label{eqn:A}
 A = S \left[ \mathcal{A} \right]_\mathcal{S} S^{-1}, 
\end{equation}
where the square matrix $S$ is formed by the vectors given in Eq.~(\ref{eqn:basis}) written as columns and
$\left[ \mathcal{A} \right]_\mathcal{S}$ is defined by Eq.~(\ref{eqn:matrix}).
The constructed 
matrix $A$ provides the prescribed norms of residual vectors at the first $q$ GMRES($m$) cycles when starting with $r_0$
and its spectrum is $\Lambda$.

\subsection{Difference with the work of Greenbaum, Pt\'{a}k, and Strako\v{s}~\cite{GreenbaumPtakStrakos:96}}

For the reader familiar with the work of Greenbaum, Pt\'{a}k, and
Strako\v{s}~\cite{GreenbaumPtakStrakos:96}, it might be tempting to obtain the
present result by pursuing the following scheme: fix $r_0$ and then
consider the first restarted GMRES cycle as the initial part of a full GMRES run
where the convergence is prescribed for the first $m$ iterations (and
arbitrarily set for the remaining $n-m$ iterations).  Then, similarly, given the starting
residual vector $r_{1}$ provided by this first cycle, construct
the next Krylov residual subspace which provides the desired convergence
following the
scheme of Greenbaum, Pt\'{a}k, and Strako\v{s}~\cite{GreenbaumPtakStrakos:96}.
Proceed
identically for the remaining cycles.  This approach, however, does not guarantee the linear independence of the
set $\overline{\mathcal{S}}$ and, hence, one meets the problem of defining the
linear operator $\mathcal{A}$. These considerations were the reason for the
assumption (\ref{eqn:simple_new_res}) on the residual reduction inside a cycle,
which allowed to quite easily justify the linear independence of the set
$\overline{\mathcal{S}}$ and, as well, to control the spectrum.

\section{Generating stagnating example of restarted GMRES}
\label{sec:stagnation}

Theorem~\ref{theo:main} handles the case for
 the decreasing positive sequence $\{f(k)\}_{k=0}^q$.
In this section, we are concerned with the stagnation case:
when $ f(0) > f(1) >  \dots  > f(s) > 0$ and $f(s) = f(s+1) =  \ldots = f(q)$.

\begin{theorem}\label{thm:1}
Given a matrix order $n$, a restart parameter $m$ ($m < n$),
a positive sequence $ \{f(k)\}_{k=0}^q$, which is either 
decreasing, or such that  
$ f(0) > f(1) >  \dots  > f(s) > 0$ and $f(s) = f(s+1) =  \ldots = f(q)$,
where $q < n/m$, $s<q$.
There exits an $n$-by-$n$ matrix $A$ and a vector $r_0$ with
$\|r_0\|=f(0)$ such that $\|r_k\|=f(k)$, $k=1,\ldots,q$, where $r_k$ is the
residual at cycle $k$ of restarted GMRES with restart parameter $m$ applied to the linear
system $Ax=b$, with initial residual $r_0=b-Ax_0$. Moreover, the matrix $A$ can be chosen
to have any desired eigenvalues.
\end{theorem}
\textbf{Proof.}
The decreasing convergence case is handled by Theorem~\ref{theo:main}.
Therefore, we only need to construct a matrix $A$ with a spectrum $\Lambda$ and an
initial residual vector $r_0$, $\|r_0\| = f(0)$ for which restarted GMRES
stagnates at cycle $s+1$ while $\|r_1\| = f(1) > \ldots > \|r_s\| = f(s)$,
$s<q$.

By Lemma~\ref{lemm:1}, given the initial residual vector $r_0$, 
$\|r_0\| = f(0)$, we can construct residual vectors $r_k$ with the prescribed norms $f(k)$, and orthonormal sets $\mathcal{W}_m^{(k)}$, $k = 1, \ldots, s$, such that the set
\begin{equation}\label{eqn:stag}        
\lbrace r_0, \mathcal{W}_{m-1}^{(1)}, \dots, r_{s-1}, \mathcal{W}_{m-1}^{(s)}, r_{s}\rbrace
\end{equation}
is linearly independent.
In order to enforce
stagnation at the ($s+1$)-st GMRES($m$) cycle, we want the next orthonormal set $\mathcal{W}_{m}^{(s+1)}$
to be orthogonal to the residual vector $r_s$.
(See Eq.~(\ref{eqn:new_res}) or~(\ref{eqn:simple_new_res}).)
Thus, following the pattern in Lemma~\ref{lemm:1}, we choose
$\mathcal{W}_{m}^{(s+1)}$ from the orthogonal complement of the set~(\ref{eqn:stag}), and append $\mathcal{W}_{m-1}^{(s+1)}$
to the set~(\ref{eqn:stag}), thus obtaining the linearly independent set     
\begin{equation}\label{eqn:stag_exp}        
\lbrace r_0, \mathcal{W}_{m-1}^{(1)}, \dots, r_{s-1}, \mathcal{W}_{m-1}^{(s)}, r_{s}, \mathcal{W}_{m-1}^{(s+1)}\rbrace. 
\end{equation}

At this point, if we followed the proof of Lemma~\ref{lemm:1}, we would
append the new residual vector $r_{s+1}$ to the
set~(\ref{eqn:stag_exp}). Since
$r_s = r_{s+1}$, this would result in the loss of the linear independence of our set.
Instead, we would like
to expand the set (\ref{eqn:stag_exp}) by some vector that will not spoil the linear
independence and will allow for a proper definition of the operator $\mathcal{A}$ at the second step of the proof.
To fulfill this task, we choose this
vector to be $w_m^{(s+1)} + r_s$ and append it to (\ref{eqn:stag_exp}). We 
obtain the set   
\begin{equation}\label{eqn:stag_exp_res}
\lbrace r_0, \mathcal{W}_{m-1}^{(1)}, \dots, r_{s-1}, \mathcal{W}_{m-1}^{(s)}, r_{s}, \mathcal{W}_{m-1}^{(s+1)},
w_m^{(s+1)} + r_s\rbrace,
\end{equation}
which is linearly independent, since the vector $w_m^{(s+1)} + r_s$ has the component $w_m^{(s+1)}$ from the orthogonal
complement of (\ref{eqn:stag_exp}) and hence cannot be represented as a linear combination of vectors in this set.

Expanding (\ref{eqn:stag_exp_res}) with vectors $\widehat{\mathcal{S}} = \lbrace
\widehat{s}_1, \dots, \widehat{s}_t \rbrace$, we finally construct the basis of $\mathbb{C}^n$:
\begin{equation}\label{eqn:basis_stag}
\tilde{\mathcal{S}} = 
\lbrace r_0, \mathcal{W}_{m-1}^{(1)}, \dots, r_{s-1}, \mathcal{W}_{m-1}^{(s)}, r_{s}, \mathcal{W}_{m-1}^{(s+1)},
w_m^{(s+1)} + r_{s}, \widehat{s}_1, \dots, \widehat{s}_t\rbrace,
\end{equation}
where $t = n - m(s + 1) - 1$.

Now, following the previously described pattern, we need to define an operator $\mathcal{A}$ with a prescribed
spectrum $\Lambda$, represented by the matrix $A$ in the canonical basis, such that 
Eq.~(\ref{eqn:span}) is satisfied for $k = 1, \dots, s+1$. We split $\Lambda$ into the disjoint subsets $\Lambda = \lbrace \Lambda_1, \Lambda_2,  \dots , \Lambda_{s+1}, \Lambda_{s+2} \rbrace$, so that each $\Lambda_k$ consists of $m$ sequential
elements of $\Lambda$, $k = 1, \dots, s+1$, while the rest $n - m(s+1)$ elements are included into $\Lambda_{s+2}$.
Similarly to (\ref{eqn:charact_m})--(\ref{eqn:charact_t}), for each $k$, we introduce the polynomials
\begin{eqnarray}
\label{eqn:charact_m_stag}   p_k(x)     & = & x^{m} - \sum_{j=0}^{m-1}\alpha_j^{(k)} x^{j}, \quad k = 1, \dots , s+1;    \\
\label{eqn:charact_t_stag}   p_{s+2}(x) & = & x^{t+1} - \sum_{j=0}^{t}\alpha_j^{(s+2)} x^{j}, \quad t = n-m(s+1)-1; 
\end{eqnarray}
where the roots of each polynomial are in the respective set $\Lambda_k$, $k = 1, \ldots, s+2$.

Similarly to (\ref{eqn:operator}), we define the operator $\mathcal{A}$ as following:   
\begin{eqnarray}
\nonumber   \mathcal{A} r_{k-1}    & = &     w_1^{(k)},       \\
\nonumber   \mathcal{A} w_1^{(k)}  & = &     w_2^{(k)},       \\
\nonumber                               &\vdots&                         \\
\nonumber   \mathcal{A} w_{m-2}^{(k)}  & = &  w_{m-1}^{(k)},                                       \\
\nonumber   \mathcal{A} w_{m-1}^{(k)}  & = &  -\alpha_0^{(k)} r_k + \alpha_0^{(k)} r_{k-1} + 
                                 \alpha_1^{(k)} w_1^{(k)} +  \dots +\alpha_{m-1}^{(k)} w_{m-1}^{(k)}, 
                                 k = 1,\dots, s;                                       \\
\nonumber   \mathcal{A} r_{s}    & = &      w_1^{(s+1)},       \\
\nonumber   \mathcal{A} w_1^{(s+1)}  & = &    w_2^{(s+1)},       \\
\label{eqn:operator_stag}                &\vdots&                         \\
\nonumber   \mathcal{A} w_{m-2}^{(s+1)}   & = &  w_{m-1}^{(s+1)},                           \\
\nonumber   \mathcal{A} w_{m-1}^{(s+1)}  & = &  - \alpha_0^{(s+1)} (w_m^{(s+1)} + r_{s}) + \alpha_0^{(s+1)} r_s + \alpha_1^{(s+1)} w_1^{(s+1)} +  \dots + \alpha_{m-1}^{(s+1)} w_{m-1}^{(s+1)},  \\
\nonumber  \mathcal{A} (w_m^{(s+1)} + r_s)  & = &  \widehat{s}_1,                                               \\ 
\nonumber  \mathcal{A} \widehat{s}_1  & = &  \widehat{s}_2,                                       \\
\nonumber                                 &\vdots&                         																		\\                  \nonumber  \mathcal{A} \widehat{s}_{t-1}  & = &  \widehat{s}_t,                                   \\
\nonumber  \mathcal{A} \widehat{s}_t  & = &  \alpha_0^{(s+2)} (w_m^{(s+1)} + r_s) + \alpha_1^{(s+2)} \widehat{s}_{1} +  \dots 
+\alpha_{t}^{(s+2)} \widehat{s}_{t},  
\end{eqnarray}                   
where $\alpha_j^{(k)}$'s are the coefficients of polynomials (\ref{eqn:charact_m_stag}) and (\ref{eqn:charact_t_stag}).
From the definition (\ref{eqn:operator_stag}) of the operator $\mathcal{A}$, one can observe that
for each $k$, $k = 1, \ldots, s+1$, 
\begin{equation*}
\mbox{span}\{\mathcal{A} r_{k-1}, \ldots, \mathcal{A}^j r_{k-1}\} = \mbox{span}~\mathcal{W}_m^{(j)}, \ j = 1, \ldots, m. 
\end{equation*}
Thus, given the representation $A$ of the operator $\mathcal{A}$ in the canonical basis, we can guarantee that Eq.~(\ref{eqn:span})
holds for each $k$, $k = 1, \ldots, s+1$. 

Similarly to Eq.~(\ref{eqn:matrix}), the
structure of the matrix $[A]_{\tilde{\mathcal{S}}}$ of the operator $\mathcal{A}$ in the basis $\tilde{\mathcal{S}}$
will be block lower triangular with
each diagonal block being the companion matrix for the corresponding subset $\Lambda_k$ of $\mathcal{A}$'s
eigenvalues, where characteristic polynomials are defined by (\ref{eqn:charact_m_stag})--(\ref{eqn:charact_t_stag}), 
and $-\alpha_0^{(k)}$'s being subdiagonal elements. The desired matrix $A$ is then obtained by similarity transformation
\begin{equation*}
 A = \tilde{S} \left[ \mathcal{A} \right]_{\tilde{\mathcal{S}}} \tilde{S}^{-1}, 
\end{equation*}
where the square matrix $\tilde{S}$ is formed by the set of vectors~(\ref{eqn:basis_stag}) written as columns.  
\begin{flushright}
 $\square$
\end{flushright}

\section{Restarted GMRES with variable restart parameter}
\label{sec:variable_restart}

The result given by Theorem~\ref{thm:1} generalizes to the case when the restart parameter $m$ is not fixed,
but varies over the successive cycles with a priori prescribed restart parameters $m_k$ for the corresponding
GMRES($m_k$) cycles.   
\begin{corollary}\label{cor:1}
Given a matrix order $n$, a sequence $\{m_k\}_{k=1}^q$ of restart parameters with $1 \leq m_k \leq n-1$,
$\displaystyle{\sum_{k=1}^q} m_k < n $, and 
a positive sequence $ \{f(k)\}_{k=0}^q$, which is either 
decreasing, or such that  
$ f(0) > f(1) >  \dots  > f(s) > 0$ and $f(s) = f(s+1) =  \ldots = f(q)$,
where $s<q$.
There exits an $n$-by-$n$ matrix $A$ and a vector $r_0$ with
$\|r_0\|=f(0)$ such that $\|r_k\|=f(k)$, $k=1,\ldots,q$, where $r_k$ is the
residual at cycle $k$ of restarted GMRES with a variable restart parameter $m_k$ applied to the linear
system $Ax=b$, with initial residual $r_0=b-Ax_0$. Moreover, the matrix $A$ can be chosen
to have any desired eigenvalues.
\end{corollary}
\textbf{Proof} The proof follows directly from Lemma~\ref{lemm:1}, Lemma~\ref{lemm:2} and Theorem~\ref{theo:main}.
Note that the constructed operator $\mathcal{A}$ will have block lower triangular matrices with block sizes $m_k$ (instead of $m$).
\begin{flushright}
 $\square$
\end{flushright}

\section{Generating non-convergent examples}
\label{sec:nonconvergent}

When constructing a matrix $A$ and an initial residual vector $r_0$ which
provide the prescribed decreasing cycle-convergence generated by GMRES($m$), we note that
from the last line of
the definition (\ref{eqn:operator}) of the operator $\mathcal{A}$
we obtain
\begin{equation*}
 r_q \in A \mathcal{K}_{t+1} \left( A, r_q\right) ,
\end{equation*}
where $A$ is the representation of the operator $\mathcal{A}$ in the canonical basis and $t = n - mq - 1$.
This equality implies that at the end of the $(q+1)$-st cycle GMRES($m$) converges to the exact solution of Eq.~(\ref{eqn:sys}), 
i.e. $r_{q+1} = 0$. This fact might seem unnatural and undesirable, e.g., for constructing academic examples.
The ``drawback'', however, can be easily fixed by a slight correction of the basis $\mathcal{S}$ -- 
somewhat similarly to how we handled the stagnation case in Theorem \ref{thm:1}.

Given residuals $r_k$ and orthonormal sets $\mathcal{W}_m^{(k)}$ constructed according
to Lemma~\ref{lemm:1},
instead of considering the set~$\mathcal{S}$,
 we consider the following basis of $\mathbb{C}^n$:
\begin{equation}\label{eqn:basis_mod}
 \mathcal{\tilde{S}} = \lbrace r_0, w_1^{(1)},  \dots , w_{m-1}^{(1)},  \ldots , r_{q-1}, 
w_1^{(q)},  \dots , w_{m-1}^{(q)}, r_q + \gamma r_{q-1}, \widehat{s}_1,  \dots , \widehat{s}_t \rbrace, 
\end{equation}
where $\gamma \neq -1$. Here we substituted the basis vector $r_q$ in Eq.~(\ref{eqn:basis}) by $r_q + \gamma r_{q-1}$. 
The vector $r_q + \gamma r_{q-1}$ cannot be represented as a linear combination of other vectors in $\tilde{S}$,
since it contains the component $r_q$, which is not represented by these vectors. Hence, $\mathcal{\tilde{S}}$ is
indeed a basis of $\mathbb{C}^n$. Thus we can define the operator
$\mathcal{A}$ by its action on $\mathcal{\tilde{S}}$:
\begin{eqnarray}
\nonumber   \mathcal{A} r_{k-1}    & = &     w_1^{(k)},       \\
\nonumber   \mathcal{A} w_1^{(k)}  & = &     w_2^{(k)},       \\
\nonumber                               &\vdots&                         \\
\nonumber   \mathcal{A} w_{m-2}^{(k)}  & = &  w_{m-1}^{(k)},                                       \\
\nonumber   \mathcal{A} w_{m-1}^{(k)}  & = &  -\alpha_0^{(k)} r_k + \alpha_0^{(k)} r_{k-1} + 
                                    \alpha_1^{(k)} w_1^{(k)} +  \dots +\alpha_{m-1}^{(k)} w_{m-1}^{(k)}, 
                                     k = 1,\dots, q-1;                                       \\
\nonumber   \mathcal{A} r_{q-1}    & = &      w_1^{(q)},       \\
\nonumber   \mathcal{A} w_1^{(q)}  & = &    w_2^{(q)},       \\
\label{eqn:operator_mod}                &\vdots&                         \\
\nonumber   \mathcal{A} w_{m-2}^{(q)}   & = &  w_{m-1}^{(q)},                           \\
\nonumber \mathcal{A} w_{m-1}^{(q)} & = & \frac{-\alpha_0^{(q)}}{1 + \gamma} (r_q + \gamma r_{q-1}) + \alpha_0^{(q)} r_{q-1} + 
                                    \alpha_1^{(q)} w_1^{(q)} + \dots 
                                    +\alpha_{m-1}^{(q)} w_{m-1}^{(q)},                                         \\
\nonumber \mathcal{A} (r_{q} + \gamma r_{q-1}) & = & \widehat{s}_1,                            \\ 
\nonumber \mathcal{A} \widehat{s}_1  & = &  \widehat{s}_2,                                       \\
\nonumber             &\vdots&                                                   \\
\nonumber \mathcal{A} \widehat{s}_{t-1}  & = & \widehat{s}_t,                                   \\
\nonumber \mathcal{A} \widehat{s}_t & = & \alpha_0^{(q+1)} (r_{q} + \gamma r_{q-1}) + \alpha_1^{(q+1)} \widehat{s}_{1} + 
           \dots  +\alpha_{t}^{(q+1)} \widehat{s}_{t},                            
\end{eqnarray}                   
where $\alpha_j^{(k)}$'s are the coefficients of the corresponding characteristic polynomials (\ref{eqn:charact_m}) and
(\ref{eqn:charact_t}). The fact that the operator $\mathcal{A}$ produces the correct Krylov residual subspace at
the cycle $q$, i.e.,
\begin{equation*}
\mbox{span}\{\mathcal{A} r_{q-1}, \dots, \mathcal{A}^m r_{q-1}\} = \mbox{span}~\mathcal{W}_m^{(q)},
\end{equation*}
can be observed from the following equalities:
\begin{eqnarray}
\nonumber \mathcal{A} w_{m-1}^{(q)} & = & \frac{-\alpha_0^{(q)}}{1 + \gamma} (r_q + \gamma r_{q-1}) + \alpha_0^{(q)} r_{q-1} + 
                                    \alpha_1^{(q)} w_1^{(q)} + \dots +\alpha_{m-1}^{(q)} w_{m-1}^{(q)}      \\
\nonumber                           & = & \frac{-\alpha_0^{(q)}}{1 + \gamma} (r_q - r_{q-1} + (1+\gamma) r_{q-1}) +                                                         \alpha_0^{(q)} r_{q-1} + \alpha_1^{(q)} w_1^{(q)} + \dots 
                                    +\alpha_{m-1}^{(q)} w_{m-1}^{(q)}                                         \\
\nonumber                           & = & \frac{-\alpha_0^{(q)}}{1 + \gamma} (r_q - r_{q-1}) + \alpha_1^{(q)} w_1^{(q)} 
                                          + \dots + \alpha_{m-1}^{(q)} w_{m-1}^{(q)},                                         
\end{eqnarray}
where, by Eq.~(\ref{eqn:operator_mod}), $\mathcal{A} w_{m-1}^{(q)} = \mathcal{A}^m r_{q-1}$ and, 
by Eq.~(\ref{eqn:simple_new_res}), $(r_q - r_{q-1}) \in \mbox{span}\{w_m^{(q)}\}$. 

The matrix $\left[ \mathcal{A} \right]_\mathcal{\tilde{S}} $ of the operator
$\mathcal{A}$, defined by Eq.~(\ref{eqn:operator_mod}), in the basis
$\mathcal{\tilde{S}}$ 
is identical to
Eq.~(\ref{eqn:matrix}) with the only change of the subdiagonal element
$-\alpha_0^{(q)}$ to $\frac{-\alpha_0^{(q)}}{1+\gamma}$, $\gamma \neq -1$.
Hence, $\mathcal{A}$ has the desired spectrum $\Lambda$.

Thus, finally, according to Eq.~(\ref{eqn:operator_mod}), 
\begin{equation*}
 r_q \in \mathcal{AK}_{t+1} \left( \mathcal{A}, r_q\right) + \mathcal{K}_{t+2} \left( \mathcal{A}, r_{q-1}\right),
\end{equation*}
providing that $r_{q+1}$ is nonzero.

\bibliography{vela09}

\begin{thebibliography}{10}

\bibitem{Baglama.Calvetti.Golub.Reichel:99}
J.~Baglama, D.~Calvetti, G.~H. Golub, and L.~Reichel.
\newblock Adaptively preconditioned {GMRES} algorithms.
\newblock {\em SIAM Journal on Scientific Computing}, 20(1):243--269, 1999.

\bibitem{BakerJessupKolev:09}
A.~H. Baker, E.~R. Jessup, and Tz.~V. Kolev.
\newblock A simple strategy for varying the restart parameter in {GMRES}(m).
\newblock {\em Journal of Computational and Applied Mathematics},
  230(2):751--761, 2009.

\bibitem{Beckermann.Goreinov.Tyrtyshnikov.05}
B.~Beckermann, S.~A. Goreinov, and E.~E. Tyrtyshnikov.
\newblock Some remarks on the {Elman} estimate for {GMRES}.
\newblock {\em SIAM Journal on Matrix Analysis and Applications},
  27(3):772--778, 2005.

\bibitem{Chapman.Saad:1997}
Andrew Chapman and Yousef Saad.
\newblock Deflated and augmented {Krylov} subspace techniques.
\newblock {\em Numerical Linear Algebra with Applications}, 4(1):43--66, 1997.

\bibitem{Duff.Giraud.Langou.Martin:05}
I.~S. Duff, L.~Giraud, J.~Langou, and {\'E}.~Martin.
\newblock Using spectral low rank preconditioners for large electromagnetic
  calculations.
\newblock {\em Int. J. Numerical Methods in Engineering}, 62(3):416--434, 2005.

\bibitem{Eisenstat.Elman.Schultz.83}
S.~C. Eisenstat, H.~C. Elman, and M.~H. Schultz.
\newblock Variational iterative methods for nonsymmetric systems of linear
  equations.
\newblock {\em SIAM Journal on Numerical Analysis}, 20:345--357, 1983.

\bibitem{Elman.82}
H.~C. Elman.
\newblock {\em Iterative methods for large sparse nonsymmetric systems of
  linear equations}.
\newblock PhD thesis, Yale University: New Haven, CT, 1982.

\bibitem{Embree.99}
M.~Embree.
\newblock How descriptive are {GMRES} convergence bounds?
\newblock Technical Report 99/08, Oxford University Computing Laboratory, 1999.

\bibitem{Embree:03}
M.~Embree.
\newblock The tortoise and the hare restart {GMRES}.
\newblock {\em SIAM Review}, 45(2):259--266, 2003.

\bibitem{Erhel.Burrage.Pohl:96}
J.~Erhel, K.~Burrage, and B.~Pohl.
\newblock Restarted {GMRES} preconditioned by deflation.
\newblock {\em Journal of Computational and Applied Mathematics}, 69:303--318,
  1996.

\bibitem{Greenbaum:97}
A.~Greenbaum.
\newblock {\em Iterative Methods for Solving Linear Systems}.
\newblock SIAM, 1997.

\bibitem{GreenbaumPtakStrakos:96}
A.~Greenbaum, V.~Pt\'{a}k, and Z.~Strako\v{s}.
\newblock Any nonincreasing convergence curve is possible for {GMRES}.
\newblock {\em SIAM Journal on Matrix Analysis and Applications},
  17(3):465--469, 1996.

\bibitem{Morgan.02}
R.~B. Morgan.
\newblock {GMRES} with deflated restarting.
\newblock {\em SIAM Journal on Scientific Computing}, 24(1):20--37, 2002.

\bibitem{NachtigalReddyTrefethen:92}
N.~M. Nachtigal, S.~C. Reddy, and L.~N. Trefethen.
\newblock How fast are nonsymmetric matrix iterations?
\newblock {\em SIAM Journal on Matrix Analysis and Applications},
  13(3):778--795, 1992.

\bibitem{SaadSchultz:86}
Y.~Saad and M.~H. Schultz.
\newblock {GMRES}: A generalized minimal residual algorithm for solving
  nonsymmetric linear systems.
\newblock {\em SIAM Journal on Scientific and Statistical Computing},
  7(3):856--869, 1986.

\bibitem{SimonciniSzyld:08}
V.~Simoncini and D.~Szyld.
\newblock New conditions for non-stagnation of minimal residual methods.
\newblock {\em Numerische Mathematik}, 109(3):477--487, 2008.

\bibitem{Simoncini.Szyld.05a}
V.~Simoncini and D.~B. Szyld.
\newblock On the occurrence of superlinear convergence of exact and inexact
  {K}rylov subspace methods.
\newblock {\em SIAM Review}, 47:247--272, 2005.

\bibitem{vanderVorstVuik:93}
H.~A. van~der Vorst and C.~Vuik.
\newblock The superlinear convergence behaviour of {GMRES}.
\newblock {\em Journal of Computational and Applied Mathematics},
  48(3):327--341, 1993.

\bibitem{VechLangou:09}
E.~Vecharynski and J.~Langou.
\newblock The cycle-convergence of restarted {GMRES} for normal matrices is
  sublinear.
\newblock SIAM Journal on Scientific Computing, to appear.

\bibitem{ZhongMorgan:08}
B.~Zhong and R.~B. Morgan.
\newblock Complementary cycles of restarted {GMRES}.
\newblock {\em Numerical Linear Algebra with Applications}, 15(6):559--571,
  2008.

\bibitem{Zitko:00}
J.~Z{\'{\i}}tko.
\newblock Generalization of convergence conditions for a restarted {GMRES}.
\newblock {\em Numerical Linear Algebra with Applications}, 7(3):117--131,
  2000.

\bibitem{Zitko:08}
J.~Z{\'{\i}}tko.
\newblock Some remarks on the restarted and augmented {GMRES} method.
\newblock {\em Electronic Transactions on Numerical Analysis}, 31:221--227,
  2008.

\end{thebibliography}
\bibliographystyle{plain}

\end{document}